\documentclass[11pt,leqno]{article}

\usepackage{amsmath,amsfonts,amscd,amssymb,theorem}

\long\def\comment#1\endcomment{}

\makeatletter
\begingroup
\gdef\th@dotted{\normalfont\itshape
  \def\@begintheorem##1##2{%
        \item[\hskip\labelsep \theorem@headerfont ##1\ ##2.]}%
\def\@opargbegintheorem##1##2##3{%
   \item[\hskip\labelsep \theorem@headerfont ##1\ ##2\ (##3).]}}
\endgroup
\makeatother

\theoremstyle{dotted}

\newtheorem{theorem}{Theorem}[section]
\newtheorem{lemma}[theorem]{Lemma}

\newtheorem{prop}[theorem]{Proposition}

\newtheorem{cond}[theorem]{Condition}

\makeatletter
\begingroup
\gdef\th@upshape{\normalfont
  \def\@begintheorem##1##2{%
        \item[\hskip\labelsep \theorem@headerfont ##1\ ##2.]}%
\def\@opargbegintheorem##1##2##3{%
   \item[\hskip\labelsep \theorem@headerfont ##1\ ##2\ (##3).]}}
\endgroup
\makeatother

\theoremstyle{upshape}

\newtheorem{defn}[theorem]{Definition}
\newtheorem{remark}[theorem]{Remark}

\makeatletter
\renewcommand{\subsection}{\@startsection{subsection}{2}{0pt}{-3ex
plus -1ex minus -0.2ex}{-2mm plus -0pt minus
-2pt}{\normalfont\bfseries}} \makeatother

\makeatletter
\@addtoreset{equation}{section}
\makeatother

\newcommand{\cntrct}                
{\hspace{2pt}\raisebox{1pt}{\text{$\lrcorner$}}\hspace{2pt}}

\newcommand{\proof}[1][Proof.]{\smallskip\noindent{\em #1}}
\def\endproof{\hfill\ensuremath{\square}\par\medskip}

\def\eqref#1{\thetag{\ref{#1}}}

\let\latexref=\ref
\def\ref#1{{\normalfont{\latexref{#1}}}}

\newcommand{\wt}{\widetilde}
\newcommand{\wh}{\widehat}

\newcommand{\ratto}{\dasharrow}     

\newcommand{\hdot}{{\:\raisebox{3pt}{\text{\circle*{1.5}}}}}
\newcommand{\C}{{\mathbb C}}
\newcommand{\Q}{{\mathbb Q}}

\newcommand{\calo}{{\mathcal O}}

\newcommand{\X}{{\mathcal X}}
\newcommand{\Y}{{\mathcal Y}}

\newcommand{\U}{{\mathcal U}}

\newcommand{\K}{{\cal K}}
\newcommand{\T}{{\cal T}}
\newcommand{\A}{{\cal A}}
\newcommand{\F}{{\cal F}}
\renewcommand{\L}{{\cal L}}
\newcommand{\N}{{\cal N}}

\newcommand{\Spec}{\operatorname{Spec}}
\newcommand{\Proj}{\operatorname{{\cal P}{\it roj}}}
\newcommand{\Pic}{\operatorname{Pic}}
\newcommand{\Ext}{\operatorname{Ext}}

\newcommand{\cchar}{\operatorname{\sf char}} 
\renewcommand{\dim}{\operatorname{\sf dim}} 
\newcommand{\codim}{\operatorname{\sf codim}}

\newcommand{\id}{\operatorname{\sf id}}
\newcommand{\Tw}{\operatorname{\sf Tw}}

\newcommand{\disc}{\operatorname{\sf disc}}
\newcommand{\gr}{\operatorname{\sf gr}}

\title{Symplectic resolutions: deformations\\ and birational maps}

\author{D. Kaledin\thanks{Partially supported by CRDF grant
RM1-2087.}}

\begin{document}

\maketitle

\begin{abstract}
We study projective birational maps of the form $\pi:X \to Y$, where
$Y$ is a normal irreducible affine algebraic variety over $\Bbb C$,
and $X$ is a smooth holomorphically symplectic resolution of the
variety $Y$. Under these assumptions, we prove several facts on the
geometry of the resolution $X$. In particular, we prove that the map
$\pi$ must be semismall. We also prove that any two such resolutions
$X_1 \to Y$, $X_2 \to Y$ are deformationally equivalent. Finally, we
prove the existence part of the log-flip conjecture for an extremal
contraction $\pi:X \to Y$ of a holomorphically symplectic projective
algebraic manifold $X$, and prove that the associated log-flip $X'$
is also symplectice and smooth.
\end{abstract}

\tableofcontents

\section*{Introduction}

In the recent papers \cite{H1}, \cite{H2}, \cite{H4} D. Huybrechts
has established some new and important facts about the geometry of
compact hyperk\"ahler manifolds. In particular, he proved the
following result:
\begin{itemize}
\item If two simple projective holomorphically symplectic manifolds
$X_1$, $X_2$ are birationally equivalent, then they are
diffeomorphic and correspond to non-separated points in the moduli
space. In other words, for such a pair $X_1$, $X_2$ there exist
one-parameter deformations $\X_1$, $\X_2$ over the base $S = \Spec
\C[[t]]$ whose generic fibers are isomorphic as complex manifolds.
\end{itemize}
Here {\em simple} means that $H^{2,0}(X_1) = \C$ is generated by the
holomorphic symplectic form $\Omega \in \Omega^2(X_1)$, and the same
is true for $X_2$.

In this paper we try to establish a {\em local} counterpart of
Huybrechts' results. Namely, we consider an affine complex variety
$X$ which admits two smooth projective resolutions $\pi_1:X_1 \to
Y$, $\pi_2:X_2 \to Y$, and we assume that the manifolds $X_1$, $X_2$
are holomorphically symplectic. In this situation, we would like to
prove that there exist one-parameter deformations $\X_1$, $\X_2$,
$\Y$ over the base $S = \Spec \C[[t]]$ and projective maps
$\pi_1:\X_1 \to \Y$, $\pi_2:\X_2 \to \Y$, such that
\begin{itemize}
\item The fibers over the special point $o \in S$ of the varieties
$\X_1/\Y$, $\X_2/\Y$ are isomorphic to the given resolutions
$X_1/Y$, $X_2/Y$.
\item The generic fibers of the varieties $\X_1$, $\X_2$ are
isomorphic.
\end{itemize}
Unfortunately (and contrary to what has been claimed in the first
version of this paper) we are only able to prove this result under a
very restrictive technical assumption (Condition~\ref{techn}).
aRoughly speaking, we need certain graded algebra to be finitely
generated. We do believe that this condition is always satisfied, so
that our result holds in full generality. However, at present we do
not have a complete proof of this fact. We plan to return to this
matter in a subsequent paper.

Assuming Condition~\ref{techn}, we will prove slightly
more. Firstly, it turns out that the generic fibers of the
deformations $\X_1$, $\X_2$ are not only isomorphic to each other,
but they are both isomorphic to the generic fiber of the deformation
$\Y$ (in particular, all these generic fibers are affine). Secondly,
we can drop the smoothness assumption on the variety $X_2$: it
suffices to require that the resolution $\pi_2:X_2 \to Y$ is crepant
(see Section~\ref{prem}, in particular Theorem~\ref{main} for
precise statements). Moreover, in some situations one can derive the
smoothness {\em a posteriori} (see Theorem~\ref{flp.smth}, proved in
Section~\ref{app}).

In the course of proving our main result, we establish some general
facts on the geometry of holomorphically symplectic manifolds $X$
equipped with a proper map $\pi:X \to Y$. We believe that these
results may be of independent interest. In particular, we prove that
such a map $\pi$ must necessarily be semismall (see
Proposition~\ref{semismall}).

We believe that our main Theorem~\ref{main} is of interest for two
reasons. In the first place, it seems that smooth crepant
resolutions of singularities which occur in complex symplectic
geometry, while non-compact, share many properties with compact
holomorphically symplectic manifolds and deserve detailed study in
their own right.  This is, for example, the point of view adopted by
A. Beauville in \cite{B}. Also, several important purely local
results have been recently proved by Y. Namikawa \cite{N1},
\cite{N2}. The local version of Huybrechts' results that we prove
here seems to further substantiate this point of view.

Secondly, while the geometric picture established in \cite{H1},
\cite{H2} never caused any doubts, there is a serious technical
mistake in \cite{H2} which invalidates much of the proof (see the
erratum in \cite{H3}). Recently D. Huybrechts was able to fix his
proof (see \cite{H4}). To do this, he used a strong new theorem of
J.-P. Demailly and M. Paun (\cite{DP}) which requires hard
analysis. We feel that the unexpected complexity of the result makes
it worthwhile to consider alternative approaches, especially those
that are more algebraic. Moreover, the combination of our methods
with those of D. Huybrechts' might yield new information.

\bigskip

We would like to note that the first version of this paper did not
contain Condition~\ref{techn} and claimed the main result in full
generality. After it was completed and posted to the e-prints server
{\tt arXiv.org}, the author had an opportunity to visit the
University of Cologne, at a kind invitation of Prof. D. Huybrechts,
and to give a long talk on his results. In the course of subsequent
discussions, several gaps in the proof were indicated to the author
by D. Huybrechts, M. Lehn and J. Wierzba. All of them but one were
minor and are fixed here; however, the remaining gap is serious and
forces us to introduce Condition~\ref{techn}. This invalidates an
application of our result to a conjecture of V. Ginzburg's which
comprised Section~\ref{app} in the first version of the paper. In
this revised version, everything about the Ginzburg's conjecture has
been removed. Everything else claimed in the first version, in
particular Proposition~\ref{semismall}, still stands. The last
Section~\ref{app} is now taken up with a new smoothness result found
by the author in the course of fixing the proofs.

The reader will find a detailed outline of the contents of the paper
and precise formulations of all the results in Section~\ref{prem}.

\subsection*{Acknowledgments.} This paper grew out of joint work with
M. Verbitsky -- his ideas and suggestions were invaluable, and he
deserves equal credit for all the results (while all the mistakes
are mine alone). I would also like to thank E. Amerik, F. Bogomolov,
O. Biquard, P. Gauduchon, V. Ginzburg, D. Huybrechts, M. Lehn,
T. Pantev, V. Shokurov and J. Wierzba for interesting and useful
discussions, both in person and by e-mail. I am particularly
grateful to J. Wierzba, D. Huybrechts and M. Lehn, who have
suggested several improvements in the proofs and helped me to notice
some gaps. Parts of the paper were prepared during the author's
stays at the Ecole Polytechnique in Paris and at the University of
Cologne. The hospitality and the stimulating atmosphere of these
institutions are gratefully acknowledged.

\section{Statements of the results.}\label{prem}

\subsection*{Notations.} Throughout the paper, we work over the field
$\C$ of complex numbers. However, apart from an occasional reference
to Hodge theory, all the proofs are algebraic, and everything works
over an arbitrary algebraically closed field $k$ of characteristic
$\cchar k = 0$. {\em Birational} will always mean dominant and
generically one-to-one. {\em Holomorphically symplectic} will in
fact mean {\em algebraically symplectic}, -- in other words,
equipped with a non-degenerate closed algebraic $2$-form. We use the
term ``holomorphically symplectic'' since it is more traditional.

\subsection*{}The main subject of this paper will be what we will call {\em
symplectic resolutions}.

\begin{defn}\label{sympl.defn}
A {\em symplectic resolution} is a pair $\langle X, Y\rangle$ of a
normal irreducible affine complex algebraic variety $Y$ and a smooth
irreducible holomorphically symplectic variety $X$ equipped with a
projective birational map $\pi:X \to Y$.
\end{defn}

Note that by the Zariski's connectedness Theorem, these assumptions
imply that
$$
Y = \Spec H^0(X,\calo_X),
$$
so that both $Y$ and the map $\pi:X \to Y$ can be recovered from the
manifold $X$. We will sometimes drop the explicit reference to the
variety $Y$ and simply call the manifold $X$ itself a symplectic
resolution. From this viewpoint, the projectivity of the map $\pi:X
\to Y$ is a condition on the manifold $X$.

\medskip

Symplectic resolutions occur naturally in many contexts -- both in
the local study of birational maps $X \to Y$ between compact
algebraic varieties and independently, as smooth resolution of
various affine singular spaces $Y$. In this regard, we note that our
notion of symplectic resolution is stronger that the notion of a
{\em symplectic singularity} recently introduced by A. Beauville
\cite{B}. The difference is the following: we require the symplectic
form $\Omega$ on $X$ to be non-degenerate everywhere, while in
Beauville's case it suffices that $\Omega$ is non-degenerate outside
of the exceptional locus of the map $\pi:X \to Y$. On the other
hand, our notion is very close to the notion of a {\em symplectic
contraction} introduced by J. Wierzba \cite{W}.

The condition of projectivity of the map $\pi:X \to Y$ is probably
too strong. We feel that most of our statements should hold in the
case when $\pi$ is proper. However, we have not been able to prove
anything in this more general case.

\medskip

We can now give an overview of the paper and formulate the main
results.

\medskip

In Section~\ref{geom}, we study geometry of a general symplectic
resolution $\pi:X \to Y$ by elementary methods. Apart from several
technical results which we need later on, we also prove here the
following statement, which might be of independent interest. The
proof is in Subsection~\ref{ss.sub}.

\begin{prop}\label{semismall}
Every symplectic resolution $\pi:X \to Y$ is semismall. In other
words, for every closed irreducible subvariety $Z \subset Y$ we have
$$
2\codim Z \geq \codim \pi(Z).
$$
\end{prop}

Section~\ref{def} presents an overview of the deformation theory for
symplectic manifolds that has been developed recently by the author
jointly with M. Verbitsky (see \cite{KV}).

In Section~\ref{tw}, we use the general deformation machinery of
Section~\ref{def} to construct one particular deformation of a
symplectic resolution $X$. This deformation $\X/S$ is over the
one-dimensional formal disc $S = \Spec\C[[t]]$, with the closed
point $o \in S$ and the generic point $\eta \in S$ (note that $\eta$
is geometrically not a point, but a punctured disc. We call the
deformation $\X/S$ {\em the twistor deformation}. Analysing its
properties, we prove the first of the main results of this paper
(for the proof see Proposition~\ref{tw.prop}).

\begin{theorem}\label{main.cor}
Let $X$ be a symplectic resolution. Then there exists flat
deformations $\X/S$, $\Y/S$ and an $S$-map $\pi:\X \to \Y$ such that
\begin{enumerate}
\item the fibers $\X_o$, $\Y_o$ over the closed point $o \subset S$
are isomorphic, respectively, to $X$ and $Y$,
\item the isomorphisms $\X_o \cong X$, $\Y_o \cong Y$ are compatible
with the maps $\pi:X \to Y$, $\X \to \Y$, and 
\item the map $\pi:\X \to \Y$ induces an isomorphism $\pi:\X_\eta
\overset{\sim}{\to} \Y_\eta$ between fibers $\X_\eta$, $\Y_\eta$
over the generic point $\eta \in S$.
\end{enumerate}
\end{theorem}

Section~\ref{crep} is the main section of the paper. Here we
partially prove a generalization of Theorem~\ref{main.cor}, which
may be thought of as a local version of D. Huybrechts' results
\cite{H2}, \cite{H1}. To formulate this statement, recall that a
birational map $\pi:X \to Y$ between normal irreducible algebraic
varieties $X$, $Y$ which admit canonical bundles $K_X$, $K_Y$ is
called {\em crepant} if the canonical map
$$
\pi^*K_Y \to K_X
$$
defined over the non-singular locus extends to an isomorphism over
the whole manifold $X$. Every symplectic resolution $\pi:X \to Y$ is
automatically crepant. Indeed, the canonical bundle $K_X$ is trivial
-- a trivialization is given by the top power of the symplectic form
$\Omega$. Since $Y$ is normal, the map $\pi:X \to Y$ is an
isomorphism over the complement $U_Y \subset Y$ to a subvariety of
$\codim \geq 2$. This implies that $K_Y$ is also trivial (in other
words, $Y$ is Gorenstein). Therefore we have $K_X \cong \pi^*K_Y$.

The main result of the paper is the following.

\begin{theorem}\label{main}
Let $\pi:X \to Y$ be a symplectic resolution.  Assume given another
normal variety $X'$ and a crepant projective map $\pi':X' \to Y$.

Assume also that the pair $X,X'/Y$ satisfies Condition~\ref{techn}.

Then there exist flat families $\X/S$, $\Y/S$, $\X'/S$ over $S =
\Spec \C[[t]]$ and projective family maps $\pi:\X \to \Y$, $\pi':\X'
\to Y$ such that
\begin{enumerate}
\item The special fibers $\X_o$, $\X'_o$, $Y_o$
are isomorphic, respectively, to $X$, $X'$, and $Y$, and these
isomorphisms are compatible with the maps $\pi$, $\pi'$.
\item The rational map 
$$
f = \pi^{-1} \circ \pi': \X' \ratto \X
$$
is an isomorphism between the generic fibers $\X_\eta$, $\X'_\eta$.
\end{enumerate}
\end{theorem}

This includes Theorem~\ref{main.cor} if we set $X' = Y$ and $\pi' =
\id$ (the identity map). 

Finally, in the last Section~\ref{app} we give an application of our
results to a problem in the Minimal Model Program for
holomorphically symplectic manifolds. To formulate it, we need to
recall some definitions. We will only give simple versions
sufficient for our purposes; the reader is advised to consult a
standard reference on the Minimal Model Program such as the paper
\cite{KMM} (in particular \cite[\S 3-2]{KMM}).

\begin{defn}\label{extr}
Let $X$ be a normal irreducible algebraic variety with trivial
canonical divisor, and let $\Delta$ be an effective $\Q$-Cartier
divisor on $X$. A birational projective map $\pi:X \to Y$ from $X$
to a normal irreducible algebraic variety $Y$ is called {\em an
extremal contraction} if
\begin{enumerate}
\item The $\Q$-Picard group $\Pic(Y) \otimes \Q$ is a $\Q$-vector
subspace of codimension $1$ in the $\Q$-Picard group $\Pic(X)
\otimes \Q$.
\item The $\Q$-Cartier divisor $\Delta$ is anti-ample with respect
to the map $\pi:X \to Y$.
\end{enumerate}
An extremal contraction $\pi:X \to Y$ is said to be {\em of flipping
type} if it is an isomoporphism outside of a closed subset $Z
\subset Z$ of codimension $\codim Z \geq 2$.

A normal irreducible algebraic variety $X'$ with trivial canonical
divisor equipped with a projective birational map $\pi':X' \to Y$ is
called the {\em flop} of the flipping-type extremal contraction
$\pi:X \to Y$ if the strict transform $\Delta' \subset X'$ of teh
divisor $\Delta \subset X$ is $\Q$-Cartier and ample with respect to
the map $\pi':X' \to Y$.
\end{defn}

It is known (\cite[Proposition 5-1-11(2)]{KMM}) that a flop is
unique, provided it exists. In fact, we must have
$$
X' = \Proj_Y\bigoplus_k \pi_*\calo(mk[\Delta])
$$
for some integer $m \geq 1$, and a flop exists if and only if the
sheaf of graded $\calo(Y)$-algebras on the right-hand side is
finitely generated. It is conjectured (\cite[Conjecture
5-1-10]{KMM}) that a flop always exists.

We can now state our last result, proved in Section~\ref{app}.

\begin{theorem}\label{flp.smth}\mbox{}
\begin{enumerate}
\item Let $X$ be a smooth quasiprojective holomorphically symplectic
manifold, let $\pi:X \to Y$ be an extremal contraction of flipping
type, and assume that the algebraic variety $Y$ is affine.  Let
$\Delta$ be an arbitrary divisor on $X$ anti-ample wioth respect to
the map $\pi:X \to Y$.

Then the one-parameter deformation $\X \to \Y$ of the contraction $X
\to Y$ provided by Theorem~\ref{main.cor} is also an extremal
contraction of flipping type, and the divisor $\Delta$ extends
canonically to a divisor on $\X$.

\item Assume in addition that both contractions $X \to Y$, $\X \to
\Y$ admit flops $\pi':X' \to Y$, $\pi':\X' \to \Y$ with respect to
the divisor $\Delta$. 

Then the contraction $\pi':\X' \to \Y$ is flat over the base of the
deformation $\X \to \Y$, and its special fiber coincides with the
contraction $X' \to Y$. Moreover, both algebraic varieties $X'$ and
$\X'$ are smooth, and the variety $X'$ is symplectic.
\end{enumerate}
\end{theorem}

\section{Geometry.}\label{geom}

\subsection{Semismallness.}\label{ss.sub} Let $\pi:X \to Y$ be a symplectic
resolution in the sense of Definition~\ref{sympl.defn}. As we have
already noted, the canonical bundles $K_X$, $K_Y$ of the varieties
$X$, $Y$ are trivial. Since $X$ is smooth, this implies that $Y$ has
canonical singularities. By \cite{E}, this in turn implies that the
singularities of the variety $Y$ are rational. However, the latter
can proved directly, without using the general theorem of \cite{E}.

\begin{lemma}\label{vnsh}
The higher direct images $R^i\pi_*\calo_X$, $i > 0$ are trivial.
\end{lemma}

\proof{} Since $K_X \cong \calo_X$, we have
\begin{equation}\label{vanish}
R^i\pi_*\calo_X = R^i\pi_*K_X = 0
\end{equation}
for $i \geq 1$ by the vanishing theorem of Grauert and
Riemenschneider \cite{GR}.
\endproof

We see that the singularities of the variety $Y$ are indeed
rational. It is well-known that this implies that the variety $Y$ is
Cohen-Macaulay, and that the first relative cohomology sheaf
$R^1\pi_*\C_X$ with coefficients in the constant sheaf $\C_X$
vanishes. Another corollary is that $H^k(X,\calo_X) = 0$ for $k \geq
1$. Since $Y$ is affine, this follows directly from \eqref{vanish}
by the Leray spectral sequence.

\medskip

We will now turn to the non-trivial geometric property of $X/Y$,
namely, Proposition~\ref{semismall}. The proof proceeds along the
standard lines used, for example, in the paper \cite{V}. Proposition
will be derived from the following fact.

\begin{lemma}\label{pb}
Let $\pi:X \to Y$ be symplectic resolution. Denote by $\Omega \in
\Omega^2(X)$ the symplectic form on the manifold $X$.  Let $\sigma:Z
\to U$ be a smooth map of smooth algebraic manifolds, and assume
given a commutative square
\begin{equation}\label{cd}
\begin{CD}
Z @>{\eta}>> X\\
@V{\sigma}VV @VV{\pi}V\\
U @>>> Y.
\end{CD}
\end{equation}
Then there exists a $2$-form $\Omega_U \in \Omega^2(U)$ on $U$ such
that 
$$
\sigma^*\Omega_U = \eta^*\Omega.
$$
\end{lemma}

Before we prove this Lemma, we would like to make the following
remark. The conclusion of the Lemma is essentially a condition on
the restriction $\eta^*\Omega$ -- namely, the claim of the Lemma
holds if and only if $\eta^*\Omega(\xi_1,\xi_2)=0$ for every two
tangent vectors $\xi_1,\xi_2$ to $Z$ at least one of which is
vertical with respect to the map $\sigma:Z \to U$. Since the map
$\sigma:Z \to U$ is smooth, this is an open condition. Therefore it
suffices to prove the claim generically on $U$. Shrinking $U$ if
necessary, we may assume that the fibered product $X \times_Y U$
admits a projective simultaneous resolution $\wh{Z}$, so that we
have a commutative diagram
$$
\begin{CD}
\wh{Z} @>>> X \times_Y U @>>> X\\
@. @VVV @VV{\pi}V\\
@. U @>>> Y.
\end{CD}
$$
Then the pair $\wh{Z} \to U$ also satisfies the assumptions of the
Lemma. Moreover, to prove the claim for an arbitrary $Z$, it
obviously suffices to consider the case $\wh{Z} = Z$. Therefore
without any loss of generality we may impose an additional assumtion
on the map $\sigma:Z \to U$:
\begin{itemize}
\item The induced map $\eta:Z \to X \times_Y U$ is generically one-to-one.
\end{itemize}
\proof[Proof of Lemma~\ref{pb}.] The short exact sequence
$$
\begin{CD}
0 @>>> \sigma^*\Omega^1(U) @>>> \Omega^1(Z) @>>> \Omega^1(Z/U) @>>>
0
\end{CD}
$$
of relative differentials for the map $\sigma:Z \to U$ induces a
three-step filtration $\sigma^*\Omega^2(U) \subset \F \subset
\Omega^2(Z)$, with $\Omega^2(Z)/\F \cong \Omega^2(Z/U)$ and
$\F/\sigma^*\Omega^2(U) \cong \sigma^*\Omega^1(U) \otimes
\Omega^1(Z/U)$. We have to prove that the $2$-form $\eta^*\Omega$ in
fact is a section of the subsheaf $\sigma^*\Omega^2(U) \subset
\Omega^2(Z)$. We will do it in two steps.

\smallskip

\proof[Step {\normalfont 1}: $\eta^*\Omega^2(U) \in \F$.]
The proof of this step is an application of an idea of J. Wierzba
\cite{W}. It suffices to prove that for every point $u \in U$, the
restriction of the form $\eta^*\Omega$ to the fiber
$$
Z_u = Z \times_U u \subset Z
$$
vanishes. To see this, consider the complex-conjugate $(0,2)$-form
$\overline{\Omega}$. Since $H^2(X,\calo_X) = 0$, the Dolbeault
cohomology class $\left[\overline{\Omega}\right]$
vanishes. Therefore the Dolbeault cohomology class
$$
\eta^*\left[\overline{\Omega}\right]\bigg|_{Z_u}
$$
also vanishes. But the manifold $Z_u$ is compact, smooth and
projective. By Hodge theory, we obtain $\eta^*\Omega = 0$ on $Z_u$. 

\smallskip

\proof[Step {\normalfont 2}: $\eta^*\Omega \in
\sigma^*\Omega^2(U)$.]  We have proved that $\eta^*\Omega$ is a
section of the sheaf $\sigma^*\Omega^1(U) \otimes
\Omega^1(Z/U)$. Therefore for every point $u \in U$ and for every
Zariski tangent vector $\xi \in T_uU$ at the point $u$ we have a
well-defined $1$-form
$$
\alpha = \eta^*\Omega \cntrct \xi
$$
on the fiber $Z_u$. We have to prove that $\alpha = 0$ for every $u
\in U$, $\xi \in T_uU$. Since $Z_u$ is a smooth projective manifold,
the $1$-form $\alpha$ is closed, and it suffices to prove that the
cohomology class $[\alpha] \in H^1(Z_u,\C)$ vanishes.

Let $X_u = X \times_Y u$ be the (possibly singular) fibered
product. By our additional assumption, the map $\eta:Z_u \to X_u
\subset X$ is generically one-to-one.

By construction, the form $\alpha$ vanishes on every non-trivial
fiber $F \subset Z_u$ of the map $\eta:Z_u \to X_u$. Indeed, for
every tangent vector $\xi'$ to the manifold $F$ we have
$$
\alpha(\xi') = \Omega(\xi,d\eta(\xi')) = 0
$$
since $d\eta(\xi') = 0$. Therefore the cohomology class $[\alpha]$
also vanishes on every fiber $F \in Z_u$. This implies that
$[\alpha] = \eta_u^*[\alpha]'$ for some class $[\alpha]' \in
H^1(X_u,\C)$. But since $Y$ has rational singularities,
$R^1\pi_*(\C) = 0$ and $H^1(X_u,\C) = 0$.
\endproof

Proposition~\ref{semismall} is deduced from Lemma~\ref{pb} by a
standard argument which we omit to save space (see \cite[Proposition
4.16]{V} or \cite[Section 3]{W}).

\subsection{Birational maps.} Let now $\pi:X \to Y$ be a symplectic
resolution, and let $\pi':X' \to Y$ be a different, possibly
singular crepant partial resolution of the same variety $Y$.
Proposition~\ref{semismall} imposes a strong constraint on the
geomtery of $X$, but it does not say anything about the geeometry of
$X'$.  Our approach is to use the geometry of $X$ to obtain results
on the geometry of $X'$. This becomes possible because of the
following lemma. It is probably a standard fact in birational
geometry; we have included a proof for the convenience of the
reader. The argument is borrowed from \cite[2.2]{H1}. The statement
is in fact more general, since it does not require $X$ to be
symplectic -- we only need it to be smooth and crepant over $Y$.

\begin{lemma}\label{emb}
Let $X$, $X'$ be two normal algebraic varieties equipped with
crepant proper birational maps $\pi:X \to Y$, $\pi':X' \to Y$ into
the same normal irreducible Gorenstein algebraic variety $Y$. Let
$$
f = (\pi')^{-1} \circ \pi:X \ratto X' \qquad 
f^{-1} = \pi^{-1} \circ \pi':X' \ratto X.
$$ 
be the natural rational maps.
\begin{enumerate}
\item The rational map $f^{-1}:X' \ratto X$ is defined and injective
on a complement $U \subset X'$ to a closed subset $Z \subset X'$ of
codimension $\codim Z \geq 2$.
\item The pullback $f^*:\Pic(X') \to \Pic(X)$ is an embedding from
the Picard group of $X'$ to the Picard group of $X$.
\item For every line bundle $L$ on $X'$ the canonical map
$$
f^*:H^0(X',L) \to H^0(X,f^*L)
$$
is an isomorphism.
\end{enumerate}
\end{lemma}

\proof{} To prove \thetag{i}, note that every birational map is {\em
defined} on the complement $U \subset X'$ to a closed subvariety of
codimension $\geq 2$; the problem is to prove that $f^{-1}:U \to X$
is an embedding. 

By assumption, the variety $X'$ is normal. Therefore, shrinking $U
\subset X'$ if necessary, we can further assume that the definition
subset $U \subset X'$ is smooth. Since the maps $X \to Y$, $X' \to
Y$ are crepant, the canonical bundles $K_X$, $K_{X'}$ are both
isomorphic to the pull-back of the canonical bundle $K_Y$.  This
means that the Jacobian of the map $f^{-1}:U \to X$ is a section of
the trivial line bundle -- in other words, a function. Denote this
function by $\Delta$.

Being a function, $\Delta$ must extend to the whole normal variety
$X'$. By the Zariski Connectedness Theorem this implies that
$\Delta$ is lifted from a function $\Delta$ on the variety $Y$.

But both maps $\pi'$ and $\pi$ are invertible over the complement
$U_Y \subset Y$ to some closed subset $X \subset Y$ of codimension
$\codim Z \geq 2$. This means that $\Delta$ is also invertible on
$U_Y$, which implies that it is invertible everywhere.

This shows that the map $f^{-1}:U \to X$ is \'etale, hence an
embedding.

\smallskip

To prove \thetag{ii}, it suffices to notice that $f \circ f^{-1}$ is
defined and identical on an open subset $U_0 \subset U$ whose
complement is of codimension $\geq 2$. Therefore $f^*$ induces an
embedding from $\Pic(U_0)$ to $\Pic(X)$. Since $X'$ is normal,
$\Pic(X')$ is a subgroup in $\Pic(U_0)$. This proves
\thetag{ii}. Moreover, we have injective maps
$(f^{-1})^*:H^0(X,f^*L) \to H^0(X',L)$, $f^*:H^0(X',L) \to
H^0(X,f^*L)$, and the composition
$$
(f^{-1})^* \circ f^* = (f^{-1} \circ f)^* = \id^*:H^0(X',L) \to
H^0(X,f^*L) \to H^0(X',L)
$$
is the identity map. This proves \thetag{iii}. 
\endproof

Lemma~\ref{emb}~\thetag{i} applies to the particular situation of a
symplectic resolution $\pi:X \to Y$ and gives a constraint on the
map $f^{-1}:X' \ratto X$. We can use it to get a constraint which
goes in the other direction.

\begin{lemma}\label{codim.3}
Let $\pi:X \to Y$ be a symplectic resolution, and let $\pi':X' \to
Y$ be a projective crepant generically one-to-one map. The
birational map
$$
f = (\pi')^{-1} \circ \pi:X \ratto X'
$$
is defined on the preimage $\pi^{-1}(U) \subset X$ of the complement
$U = Y \setminus Z \subset Y$ to a closed subset $Z \subset Y$ of
codimension $\codim Z \geq 3$.
\end{lemma}

\proof{} Let $W \subset X$ be the indeterminacy locus of the map
$f:X \ratto X'$, and let $Z = \pi(W) \subset Y$. To prove the Lemma,
it suffices to prove that $\codim Z \geq 3$. Since $\codim W \geq
2$, the only possibility that we have to exclude is the following:
\begin{itemize}
\item There exists an irreducible component $W_0 \subset W \subset
X$ of codimension $\codim W_0 = 2$ which is generically finite over
its image $\pi(W_0) \subset Y$.
\end{itemize}
Let $W_0 \subset W \subset X$ be such a component. Consider the
graph $\Gamma \subset X \times_Y X'$ of the birational map $f:X
\ratto X'$ and let $\pi_1:\Gamma \to X$, $\pi_2:\Gamma \to X_2$ be
the canonical projections. Since $W_0$ lies in the indeterminacy
locus $W \subset X$, its preimage $\pi_1^{-1}(W) \subset \Gamma$ is
of dimension strictly greater than $\dim W_0$. Then $\codim W_0 = 2$
implies that $\pi_1^{-1}(W)\subset\Gamma$ is a divisor. 

Denote by $E = \pi_2(\pi_1^{-1}(W_0)) \subset X'$ the image of the
divisor $\pi_1^{-1}(W_0) \subset \Gamma$ under the second projection
$\pi_2:\Gamma \to X'$. By definition $\pi_1^{-1}(W)$ lies in the
fibered product $W_0 \times_Y X'$. Therefore the fibers of the map
$\pi_2:\pi^{-1}(W_0) \to E \subset X'$ are closed subschemes of the
fibers of the map $\pi:W_0 \to \pi(W_0) \subset Y$. Since $W_0$ is
generically finite over $\pi(W_0)$, this implies that
$\pi_1^{-1}(W_0)$ is generically finite over $E$.  Thus $E \subset
X'$ is also a divisor. This implies that the rational map $X' \ratto
X$ is defined in the generic point of $E$, -- in other words, that
the projection $\pi_2:\pi_1^{-1}(W_0) \to E$ is generically
one-to-one.

But by Lemma~\ref{emb}~\thetag{i} the rational map $X' \to X$ is not
only defined in the generic point of the divisor $E \subset X'$, but
also injective in this generic point. Therefore the projection
$\pi_1:\pi^{-1}_1(W_0) \to W_0$ is generically injective, which is a
contradiction. 
\endproof

\section{Deformation theory.}\label{def}

In order to study symplectic resolutions, we need an appropriate
version of Kodaira-Spencer deformation theory. We will use the
theory proposed recently by the author jointly with M. Verbitsky
\cite{KV}. We will not use the results of \cite{KV} in full
generality, but only in a very particular case. To make the present
paper self-contained, we have decided to state here all the
necessary facts, while omitting proofs and unnecessary details.  In
order to make the material more accessible by placing it in a
familiar context, we have also included a reminder on some known
facts from deformation theory, in particular on some results of
Z. Ran. The reader will find all the omitted details and much more
in the paper \cite{KV}.

\subsection*{Notation}
For every integer $n \geq 0$, denote by $S_n = \Spec\C[t]/(t^{n+1})$
the $n$-th order infinitesimal neighborhood of the closed point $o
\in S$ of the formal disk $S = \Spec\C[[t]]$. A {\em deformation of
order $n$} of a complex manifold $X$ is by definition a flat variety
$X_n/S_n$ whose fiber over $o \in S_n$ is identified with $X$.
Since $S_k \subset S_n$ for $k \leq n$, every deformation $X_n/S_n$
of order $n$ defines by restriction a deformation $X_k/S_k$ of order
$k$ for every $k \leq n$. By a {\em compatible system} of
deformations we will understand a collection of deformations
$X_n/S_n$, one for each $n \geq 1$, such that $X_n/S_n$ restricted
to order $k \leq n$ coincides with the given $X_k/S_k$.

\subsection*{}
In \cite{Bogo}, F. Bogomolov has proved the following remarkable
fact:
\begin{itemize}
\item Let $X$ be a smooth projective holomorphically symplectic
complex manifold. Then every first-order deformation $X_1/S_1$ of
$X$ extends to a compatible system of deformations $X_n$ of all
orders $n \geq 1$.
\end{itemize}
This fact has been generalized to Calabi-Yau manifolds by G. Tian
\cite{Tian} and A. Todorov \cite{Todo}. Later, a purely algebraic
proof has been given by Z. Ran.

We will now describe this proof. Recall that the general deformation
theory associates to an $n$-th order deformation $X_n/S_n$ a certain
class
$$
\theta_n \in H^1(X_n, \T(X_n/S_n) \otimes \Omega^1(S_n)),
$$ 
called {\em the Kodaira-Spencer class} (here $\T(X_n/S_n)$ is the
relative tangent bundle, and $\Omega^1(S_n)$ is the
$\calo(S_n)$-module of K\"ahler differentials over $\C$). 

The class $\theta_n$ measures the non-triviality of the deformation
$X_n/S_n$. It would be natural to try describe all the deformations
$X_n/S_n$ in terms of the associated Kodaira-Spencer classes
$\theta_n$. However, this is not possible, since the cohomology
group that contains $\theta_n$ itself depends on the deformation
$X_n/S_n$.

The key observation in Ran's approach is that this difficulty can be
circumvented if one studies deformations step-by-step, by going
through the embeddings $i_n:S_n \hookrightarrow S_{n+1}$. Indeed,
assume given a deformation $X_n/S_n$ of order $n$, and assume that
the deformation $X_n$ is further extended to a deformation
$X_{n+1}/S_{n+1}$. Then we can restrict the Kodaira-Spencer class
$\theta_{n+1}$ to $S_n$ and obtain a cohomology class
$$
\wh{\theta}_n = i_n^*\theta_{n+1} \in H^1(X_n,\T(X_n/S_n) \otimes
i_n^*\Omega^1(S_{n+1})).
$$
Applying the canonical projection $i_n^*\Omega^1(S_{n+1}) \to
\Omega^1(S_n)$ to this class $\wh{\theta}_n$ gives the order-$n$
Kodaira-Spencer class $\theta_n$. However, the $\calo(S_n)$-module
$i_n^*\Omega^1(S_{n+1})$ is strictly bigger than the module
$\Omega^1(S_n)$ -- we have a short exact sequence
\begin{equation}\label{diff.cd}
\begin{CD}
0 @>>> \C \cdot d(t^{n+1}) @>>> i_n^*\Omega^1(S_{n+1}) @>>>
\Omega^1(S_n) @>>> 0.
\end{CD}
\end{equation}
Thus on one hand, the class $\wh{\theta}_n$ lies in a cohomology
group defined entirely in terms of $X_n$, while on the other hand,
this class carries information about the next order deformation
$X_{n+1}$.

Formalizing this, Ran proves the following (where ``lifting'' is
taken with respect to the canonical surjection in \eqref{diff.cd}).

\begin{lemma}[Ran's criterion]\label{ran}
Assume given an $n$-th order deformation $X_n/S_n$ of
the complex manifold $X$. Denote by
$$
\theta_n \in H^1(X_n, \T(X_n/S_n) \otimes \Omega^1(S_n))
$$
the associated Kodaira-Spencer class.

Then the deformation $X_n/S_n$ extends to a deformation
$X_{n+1}/S_{n+1}$ of order $n+1$ if and only if the class $\theta_n$
lifts to a class
$$
\wh{\theta}_n \in H^1(X_n,\T(X_n/S_n) \otimes
i_n^*\Omega^1(S_{n+1})).
$$
Moreover, such liftings are in one-to-one correspondence with the
isomorphism classes of the extended deformations
$X_{n+1}/S_{n+1}$.\endproof
\end{lemma}

\begin{remark}
Note that for $n=1$ this Lemma reduces to the standard fact:
deformations of order $1$ are classified, up to an isomorphism, by
cohomology classes in the group
$$
H^1(X_0,\T(X_0/S_0) \otimes i_0^*\Omega^1(S_1)) = H^1(X,\T(X)).
$$
\end{remark}

Lemma~\ref{ran} corresponds to what Ran called {\em $T_1$-lifting
property}. Having established this property, Ran then uses Hodge
theory to construct all the necessary liftings. This proves the
theorem.

In the paper \cite{KV}, this proof of the Bogolomov Theorem has been
generalized to a certain class of non-compact holomorphically
symplectic manifolds. We will not need the full generality of
\cite{KV}. It suffices to know that the results apply to
holomorphically symplectic manifolds $X$ with $H^k(X,\calo(X)) = 0$
for $k \geq 1$. By Lemma~\ref{vnsh}, this holds for all symplectic
resolutions $\pi:X \to Y$.

To pay for such a generalization, one has to change the notion of a
deformation. Namely, introduce the following.

\begin{defn}
A {\em symplectic deformation} $\X/S$ of a holomorphically
symplectic complex algebraic manifold $\langle X, \Omega \rangle$
over a local Artin base $S$ is a smooth family $\X/S$ equipped with
a closed non-degenerate relative $2$-form
$$
\Omega \in \Omega^2(\X/S),
$$
such that the fiber over the closed point $o \in S$ is identified
with $\langle X,\Omega \rangle$.
\end{defn}

In the theory of symplectic deformations, the tangent sheaf
$\T(\X/S)$ is replaced by a certain complex
$\K^\hdot(\X/S)$ of sheaves on $\X$. This complex
$$
\K^\hdot(\X/S) = F^1(\Omega^\hdot(\X/S))[-1]
$$
is the first term of the Hodge (=stupid) filtration of the relative
de Rham complex $\Omega^\hdot(\X/S)$, shifted by $-1$, -- in
other words, we set
$$
K^i(\X/S) = 
\begin{cases} \Omega^{i+1}(\X/S), &\quad i \geq 0,\\
0, &\quad i < 0,
\end{cases}
$$
with the de Rham differential. The cohomology of the complex
$\K^\hdot(\X/S)$ can be included into a long exact sequence
$$
\begin{CD}
H_{DR}^k(\X/S) @>>> H^k(\X,\calo(\X)) @>>> H^k(\X,\K^\hdot(\X/S))
@>>> \dots
\end{CD}
$$
where $H_{DR}^k(\X/S)$ is the relative de Rham cohomology of $\X$
over $S$. By assumption the coherent cohomology $H^k(\X,\calo(\X))$
vanish for $k \geq 1$. Therefore for $k \geq 1$ we have isomorphisms
$$
H^k(\X,\K^\hdot(\X/S)) \cong H_{DR}^{k+1}(\X/S).
$$
The Gauss-Manin connection $\nabla$ thus induces a connection on the
higher cohomology groups of the complex $\K^\hdot(\X/S)$.

The role of the Kodaira-Spencer class in the symplectic theory is
played by the class
$$
\theta = \nabla(\Omega) \in H^1(\X,\K^\hdot(\X/S))
\otimes \Omega(S)
$$
obtained by applying the Gauss-Manin connection to the class of the
relative symplectic form $\Omega \in \Omega^2(\X/S)$. In
particular, for an order-$n$ symplectic deformation $X_n/S_n$ we
obtain a canonical class
$$
\theta_n \in H^1(X_n,\K^\hdot(X_n/S_n) \otimes \Omega^1(S_n)).
$$
In the sequel, the only results of \cite{KV} that we will need are
the following two properties of the class $\theta_n$.

\begin{lemma}\label{ran.s}\mbox{}
\begin{enumerate}
\item The image $\tau(\theta_n) \in
H^1(X_n,\T(X_n/S_n)\otimes\Omega^1(S_n))$ of the class $\theta_n$
under the canonical map
$$
\tau:\K^\hdot(X_n/S_n) \to \K^0(X_n/S_n) \cong \Omega^1(X_n/S_n)
\cong \T(X_n/S_n))
$$
is the usual Kodaira-Spencer class for the deformation $X_n/S_n$.
\item The Ran's Criterion (Lemma~\ref{ran}) holds literally for
symplectic deformations, with the relative tangent bundle
$\T(X_k/S_k)$ replaced with the complex
$\K^\hdot(X_k/S_k)$.\endproof
\end{enumerate}
\end{lemma}

\section{Twistor deformations.}\label{tw}

\subsection{The construction.} Let $X$ be a smooth projective
symplectic resolution. To construct symplectic deformations of the
manifold $X$, we need a supply of cohomology classes in the group
$$
H^1(X,\K^\hdot(X)).
$$
Such a supply is provided by Chern classes of line bundles on
$X$. Indeed, recall that the Chern class $c_1(L) \in H^{1,1}(X)$ of a
line bundle $L$ on $X$ equals
$$
c_1(L) = d\log([L]),
$$
where $[L] \in H^1(X,\calo_X^*) = \Pic(X)$ is the class of $L$ in
the Picard group. But the $d\log$ map obviously factors through a
map
$$
d\log:\calo_X^* \to \K^\hdot(X) \to \K^0(X) = \Omega^1(X).
$$
Therefore for every line bundle $L$ we have a canonical Chern class
$$
c_1(L) \in H^1(X,\K^\hdot(X))
$$
which projects to the usual Chern class under the map
$\tau:\K^\hdot(X) \to \Omega^1(X)$. We will call it the {\em
symplectic Chern class} $c_1(L)$ of the line bundle $L$. More
generally, if we are given a deformation $X_n/S_n$ of order $n$ and
a line bundle $L$ on $X_n$, then this construction gives a canonical
class
$$
c_1(L) \in H^1(X_n,\K^\hdot(X_n/S_n))
$$
in relative cohomology.

\begin{defn}\label{tw.def}
Let $\pi:X \to Y$ be a symplectic resolution, and let $L$ be an
ample line bundles on $X$. A {\em twistor deformation} associated to
the line bundle $L$ is a compatible system of symplectic
deformations $X_n/S_n$ equipped with line bundles $L$ such that
\begin{enumerate}
\item The fiber over $o \in S_n$ of the system $\langle X_n, L
\rangle$ is identified with the pair $\langle X,L \rangle$.
\item The Kodaira-Spencer class
$$
\theta_n \in H^1(X_n,\K^\hdot(X_n/S_n) \otimes \Omega^1(S_n))
$$
of the symplectic deformation $X_n/S_n$ is equal to $c_1(L) \otimes
dt$. 
\end{enumerate}
\end{defn}

This notion is useful because of the following.

\begin{lemma}\label{pic.rigid}
Let $X$ be a symplectic resolution. For every deformation $X_n/S_n$
of an arbitrary order $n$, every line bundle $L$ on $X$ extends to a
line bundle $L$ on $X_n$, and this extension is unique up to an
isomorphism.
\end{lemma}

\proof[{Sketch of the proof (compare \cite[2.3]{H1} and
\cite[Corollary 4.3]{H1}.)}] The deformation theory of the pair
$\langle X, L \rangle$ is controlled by the (first and second)
cohomology of the Atiyah algebra $\A(L)$ of the line bundle $L$,
while the deformation theory of the manifold $X$ itself is
controlled by the (first and second) cohomology of the tangent
bundle $\T(X)$. We have a short exact sequence
$$
\begin{CD}
0 @>>> \calo_X @>>> \A(L) @>>> \T(X) @>>> 0,
\end{CD}
$$
which induces a long exact sequence of the cohomology groups.  But
since $X$ is a symplectic resolution, we have $H^i(X,\calo_X) = 0$
for $i > 0$. Therefore cohomology groups $H^i(X,\A(L))$ and
$H^i(X,\T(X))$ are isomorphic
$$
H^i(X,\A(L)) \cong H^i(X,\T(X))
$$
in every degree $i \geq 1$. \endproof

This lemma shows that the (isomorphism classes of the) line bundles
$L$ on $X_n$ are completely defined by their restrictions to the
special fiber $X \subset X_n$. Combining it with
Lemma~\ref{ran.s}~\thetag{ii}, we immediately get the following.

\begin{lemma}
For every symplectic resolution $X$ and an ample line bundle $L$ on
$X$, there exists a twistor deformation $\X/S$ associated to the
pair $\langle X,L \rangle$. Moreover, such a twistor deformation is
unique up to an isomorphism.
\end{lemma}

\proof{} Lemma~\ref{ran.s}~\thetag{ii} and induction on $n$.
\endproof

\subsection{Algebraization.} Taking the direct limit
$$
X_\infty = \lim_{\rightarrow} X_n
$$
of the order-$n$ twistor deformations $X_n/S_n$, we obtain a formal
scheme $X_\infty$ over the power series algebra $\C[[t]]$. This
formal scheme is smooth over $\C[[t]]$. It carries a line bundle $L$
and a relative symplectic form $\Omega \in \Omega^2(X_\infty/S)$.

We can now use the ampleness condition on $L$ to obtain not just the
formal scheme, but an actual deformation over the formal disc $S =
\Spec\C[[t]]$.

\begin{prop}\label{alg}
There exists a Noetherian affine scheme $\Y/S$ which is normal,
irreducible and flat over the base $S$, a smooth scheme $\X/S$ also
flat over $S$, a projective birational $S$-scheme map $\pi: \X \to
\Y$ of finite type, a line bundle $\L$ on $\X$ very ample with
respect to the map $\pi$ and a relative symplectic form $\Omega \in
\Omega^2(\X/S)$ such that the formal scheme $\langle X_\infty, L
\rangle$ with the bundle $L$ is the completion of the pair $\langle
\X,\L\rangle$ near the fiber over the special point $o \in S$, and
the form $\Omega$ on $X_\infty$ is the completion of the form
$\Omega \in \Omega^2(\X/S)$.
\end{prop}

\proof{} For every $n \geq 0$, the embedding $X_n \subset X_{n+1}$
induces a short exact sequence
$$
\begin{CD}
0 @>>> t^{n+1}\calo(X) @>>> \calo(X_{n+1}) @>>> \calo(X_n) @>>> 0
\end{CD}
$$
of sheaves on $X$, which in turn induces a cohomology long exact
sequence
$$
\begin{CD}
H^0(X,\calo(X_{n+1})) @>>> H^0(X,\calo(X_n)) @>>> H^1(X,\calo(X))
@>>>
\end{CD}
$$
Since $H^1(X,\calo(X)) = 0$ by Lemma~\ref{vnsh}, the canonical map
$$
H^0(X,\calo(X_{n+1})) \to H^0(X,\calo(X_n))
$$ 
is surjective for every $n \geq 0$. Consider the $\C[[t]]$-algebra
$$
A = \lim_{\leftarrow} H^0(X_n,\calo(X_n))
$$
and let $\Y = \Spec A$. The algebra $A$ is complete with respect to
the $(t)$-adic topology. The associated graded quotient $\gr A$ with
respect to the filtration by powers of the ideal $t \cdot A \subset
A$ is isomorphic
$$
\gr A \cong \calo(Y)[t]
$$
to the algerba of polynomials over the algebra $\calo(Y)$ of
functions on the affine variety $Y$. Therefore the algebra $A$ is
Noetherian (\cite[0, Corollaire 7.2.6]{EGA}). Moreover, $A$
obviously has no $t$-torsion, so that $\Y$ is flat over $S =
\Spec\C[[t]]$.

Let $Y_\infty$ be the completion of $\Y$ near the special fiber $Y
\subset \Y$. Then the formal scheme $X_\infty$ is proper over the
formal scheme $Y_\infty$, and the line bundle $L$ extends to a
bundle on the formal scheme $\X$. Therefore we can apply the
Grothendieck Algebraization Theorem \cite[III, Th\'eor\`eme
5.4.5]{EGA} which claims the existence of the scheme $\X/\Y$ and the
bundle $\L$ satisfying the conditions of the Proposition, except
possibly for the smoothness of $\X$ and the normality of
$\Y$. Moreover, we have
$$
\Y = \Spec H^0(\X,\calo(X)),
$$
which implies that $\pi:\X \to \Y$ is birational and that
$\pi_*\calo(\X) \cong \calo(\Y)$. 

It remains to prove that $\X$ is smooth and that $\Y$ is
normal. Since $\pi_*\calo(\X) \cong \calo(\Y)$, the second follows
from the first. Thus we have to prove that the scheme $\X$ is
regular at every point $x \in \X$. It is well known that regularity
is stable under generization. Therefore it suffices to consider only
closed points $x \in \X$. 

Fix such a point $x \in \X$. If the point $x$ lies in the special
fiber $X \subset \X$, then the regularity on $\X$ at $x$ is
equivalent to the regularity of the formal scheme $X_\infty$ at $x$.
This immediately follows from the construction of $X_\infty$. If the
point $x$ lies in the generic fiber $\X_\eta \subset \X$, then it
has no direct relation to the formal scheme, and {\em a priori}
there is no reason for $\X$ to be regular at $x \in \X$.

However, the latter situation cannot occur. Indeed, since the map
$\pi:\X \to \Y$ is proper, the point $\pi(x) \subset \Y$ is also
closed. Therefore the function $t$ either vanishes at the point $x$,
or it maps to an invertible element in the residue field
$k(\pi(x))$. In the latter case, the ideal $t\calo_\Y \subset \Y$
maps to the whole $k(\pi(x))$ under the projection $\calo_\Y \to
k(\pi(x))$, so that we have $1 = at \in k(\pi(x))$ for some function
$a \subset \calo_\Y$. Then $1-at \in \calo_\Y$ vanishes at $x$.  But
since the ring $\calo_\Y$ is complete with respect to the $(t)$-adic
topology, the function $1-at$ is invertible in $\calo_\Y$ for every
$a \in \calo_\Y$. This is a contradiction. We conclude that the
function $t$ vanishes at every closed point $x \in \X$, which means
that the point $x \in \X$ lies in the special fiber $X \subset \X$.

Finally, the existence of the form $\Omega \in \Omega^2(\X/S)$
satisfying the conditions of the Proposition is immediate from
\cite[III, Proposition 5.1.2]{EGA}.\endproof

Note that this statement is not a mere formality -- on the contrary,
it has direct geometric meaning. Indeed, while the formal scheme
$X_\infty$ is concentrated set-theoretically on the special fiber $X
\subset X_\infty$, the scheme $\X/S$ has a perfectly well-defined
generic fiber over the generic point $\eta \in S$. To pay for this,
we have to introduce the scheme $\Y$ which is {\em not of finite
type} over $S$. This has certain counterintuitive corollaries. For
instance, the construction of the twistor deformation is essentially
local on $Y$ -- if we have an open affine subset $Y_0 \subset Y$,
then the pullback $X_0 = X \times_Y Y_0$ admits a twistor
deformation $\X_0/\Y_0$, and we have canonical comparison maps
$$
\X_0 \to \X, \qquad\qquad\qquad \Y_0 \to \Y.
$$
However, it is unlikely that these maps are open embeddings.Thus one
cannot hope to obtain twistor deformations over non-affine varieties
$Y$ by working locally on $Y$ and then glueing the pieces together.

\begin{remark}
The name {\em twistor deformation} comes from hyperk\"ahler
geometry: if the holomorphically symplectic manifold $X$ admits a
compatible hyperk\"ahler metric, then the twistor deformation
$X_\infty$ is simply the completion of the twistor space $\Tw(X)$
along the fiber over the point $0 \in \C P^1$. It might be possible
to reverse the construction and obtain the whole twistor space
$\Tw(X)$ from the twistor deformation $X_\infty$, thus giving a
purely algebraic construction of a hyperk\"ahler metric on the
non-compact manifold $X$. Proposition~\ref{alg} can be considered as
the first step in this direction. However, at present it is unclear
how to proceed any further.
\end{remark}

\subsection{The generic fiber.} We can now prove the following
Proposition, which immediately implies Theorem~\ref{main.cor}.

\begin{prop}\label{tw.prop}
Let $\pi:X \to Y$ be a symplectic resolution equipped with an ample
line bundle $L$. Consider the twistor deformation $X_\infty$
associated to the line bundle $L$, and let $\pi:\X \to \Y$ be the
family of schemes over $S = \Spec\C[[t]]$ constructed from
$X_\infty$ in Proposition~\ref{alg}.

Then the map $\pi$ induces an isomorphism
$$
\pi:\X_\eta \to \Y_\eta
$$
between the fibers over the generic point $\eta \in S$.
\end{prop}

\proof{}\footnote{The argument uses a standard trick which probably
goes back to \cite{F}. Compare \cite[Propositon 4.1]{H1} and
\cite[Section 2, Claim 3 on p. 18]{N1}.}  The generic fiber $\Y_\eta
= \Y \times_S \eta$ is open in the normal scheme $\Y$ (indeed, it is
the complement to the closed special fiber $Y \subset \Y$. Since the
scheme $\Y$ is normal, the open subscheme $\Y_\eta \subset \Y$ is
also normal.  Since the map $\pi:\X \to \Y$ is birational, it
suffices to prove that it is finite over the generic point $\eta \in
S$. Moreover, by construction the map $\pi:\X_\eta \to \Y_\eta$ is
projective. Therefore by \cite[IV, Th\'eor\`eme 8.11.1]{EGA} it
suffices to prove that it is quasifinite -- in other words, that its
fibers do not contain any compact curves.

Let $\iota:C \to \X_\eta$ be an arbitrary map from a compact curve
$C/\eta$ to $\X_\eta$. Replacing $C$ with its normalization, we can
assume that the curve $C/\eta$ is connected and smooth. Since the
point $\eta = \Spec \C((t))$ is in fact the punctured disk, the
curve (or rather, the family of curves) $C/\eta$ has a well-defined
Kodaira-Spencer class
$$
\theta_C \in H^1(C,\T(C/\eta) \otimes \Omega^1(\eta/\C)) \cong
H^1(C,\T(C/\eta)).
$$
Moreover, by functoriality of the Kodaira-Spencer classes, for every
$k$-form $\alpha \in \Omega^k(\X/S)$ we have
$$
\iota^*\alpha \cntrct \theta_C = \iota^*(\alpha \cntrct
\theta_\infty) \in H^1(C,\Omega^{k-1}(C/\eta)),
$$
where $\theta_\infty \in H^1(\X,\T(X_\infty/S))$ is the
Kodaira-Spencer class of the family $\X/S$. In particular, we have
$$
\iota^*(\Omega \cntrct \theta_\infty) = \iota^*\Omega \cntrct
\theta_C,
$$
where $\Omega$ is the symplectic form. Since $C$ is a curve, we have
$\iota^*\Omega = 0$. We deduce that $\iota^*(\Omega \cntrct
\theta_\infty) = 0$. But by definition of the twistor deformation we
have
$$
\Omega \cntrct \theta_\infty = \tau(\theta_\infty) = c_1(L).
$$
Therefore $\iota^*c_1(L) = 0$. Since the line bundle $L$ is ample,
this is possible only if $\iota:C \to X_\eta$ is a projection onto a
point in the generic fiber $X_\eta$.
\endproof

\section{The main theorem.}\label{crep}

We now turn to the proof of our main Theorem~\ref{main}.

\subsection{The setup and the construction.}\label{stmnt} 
The setup for proving Theorem~\ref{main} is as follows. Assume given
a projective symplectic resolution $\pi:X \to Y$ and a different
irreducible normal variety $X'$ equipped with a projective
birational crepant map $\pi':X' \to Y$. Fix a relatively very ample
line bundle $L$ on $X/Y$ and a relatively very ample line bundle
$L'$ on $X'/Y$, so that we have
$$
X = \Proj_Y\bigoplus_{k \geq 0}\pi_*L^{\otimes k},\qquad\qquad
X' = \Proj_Y\bigoplus_{k \geq 0}\pi'_*(L')^{\otimes k},
$$
where $\Proj$ is taken on the variety $Y$.

By virtue of Lemma~\ref{emb}~\thetag{ii}, the Picard group
$\Pic(X')$ is canonically embedded into the Picard group
$\Pic(X)$. Therefore every line bundle on $X'$ defines a line bundle
on $X$. By abuse of notation, we will denote both by the same
letter. Then by Lemma~\ref{emb}~\thetag{iii} we have
\begin{equation}\label{x.prime}
X' = \Proj_Y\bigoplus_{k \geq 0}\pi'_*(L')^{\otimes k} =
\Proj_Y\bigoplus_{k \geq 0}\pi_*(L')^{\otimes k}.
\end{equation}
In other words, global sections of the bundles $(L')^{\otimes k}$
are the same, whether computed on $X$ or on $X'$.

We denote by $\pi:\X \to \Y$ the twistor deformation of the pair
$\langle X, L\rangle$, whose existence is guaranteed by
Proposition~\ref{alg}. We know that $\Y$ is a Noetherian affine
scheme, flat over $S = \Spec\C[[t]]$, $\X$ is smooth over $S$ and
projective over $\Y$, and the fiber of the map $\pi:\X \to \Y$ over
the special point $o \in S$ is identified with the given symplectic
resolution $\pi:X \to Y$. Moreover, the line bundle $L$ extends to a
relative very ample line bundle $\L$ on $\X$, so that we have
$$
\X = \Proj_{\Y}\bigoplus_{k \geq 0}\pi_*\L^{\otimes k}.
$$
The fiber $\pi:\X_\eta \to \Y_\eta$ of the map $\pi:\X \to \Y$ over
the generic point $\eta \in S$ is an isomorphism.

We can now introduce the technical condition that we need in order
to prove Theorem~\ref{main}. By Lemma~\ref{pic.rigid}, the line
bundle $L'$ on $X$ canonically extends to a line bundle $\L'$ on the
twistor deformation $\X$.

\begin{cond}\label{techn}
Consider the line bundle $\L'$ on the twistor deformation $\X$. The
graded algebra
$$
\bigoplus_{k \geq 0} H^0(\X,{\L'}^{\otimes k})
$$
is finitely generated over the algebra
$H^0(\X,\calo_{\X})=\calo(\Y)$.
\end{cond}

This will be our standing assumption from now on and till the end of
this Section.

To prove Theorem~\ref{main} under Condition~\ref{techn}, we first
construct some partial resolution $\X' \to \Y$, then prove that it
is indeed the resolution that we need. The construction is as
follows. Consider the scheme
$$
\X'_0 = \Proj_{\Y}\bigoplus_{k \geq 0} H^0(\X',{\L'}^{\otimes k})
$$
In other words, let $\X'_0$ be the image of the (rational)
projective embedding of the scheme $\X$ defined by the line bundle
$\L'$. Notice that since $\pi:\X \to \Y$ is an isomorphism over the
generic point $\eta \in S$, the scheme $\X'_0/\Y$ is isomorphic to
$\Y$ over the generic point $\eta \in S$.

\begin{lemma}\label{flt}
The scheme $\X'_0$ is reduced and flat over $S$.
\end{lemma}

\proof{} Since the scheme $\X$ is irreducible, the $\calo(\Y)$-algebra
$$
\bigoplus_{k \geq 0}\pi_*{\L'}^{\otimes k}
$$
has no zero divisors. This proves that $\X'_0$ is reduced. Moreover,
this proves that the structure sheaf $\calo(\X'_0)$ has no
$t$-torsion, which implies that $\X'_0$ is flat over $S = \C[[t]]$.
\endproof

\begin{lemma}\label{irr}
The fiber $X'_{0,o} \subset \X'_0$ of the scheme $\X'_0/S$ over the
special point $o \in S$ is reduced and irreducible.
\end{lemma}

\proof{} Since $Y$ is affine, it suffices to prove that the graded
algebra
$$
\bigoplus_{k \geq 0}\pi_*(\L')^{\otimes k}/t
$$ 
has no zero divisors. By the base change theorem for the flat map
$\X \to S$, the canonical map
$$
\bigoplus_{k \geq 0}H^0\left(\X,(\L')^{\otimes k}\right)/t \to
\bigoplus_{k \geq 0}H^0\left(X,(L')^{\otimes k}\right)
$$
is injective. Therefore it suffices to prove that the algebra on the
right-hand side has no zero divisors. But this is immediate, since
the scheme $X$ is reduced and irreducible.  
\endproof

We now let $\X'$ be the normalization of the scheme $\X'_0$, and
denote by $X'_o \subset \X'$ the special fiber of the scheme
$\X'/S$.

\medskip

By construction, $\X'$ is a normal reduced irreducible scheme
equipped with a canonical projective birational map $\pi':\X' \to
\Y$. This projection factors through a finite generically one-to-one
map $\nu:\X' \to \X'_0$. By abuse of notation, we will denote by
$$
\L' = \nu^*\calo(1)
$$
the preimage of the relative $\calo(1)$-sheaf on $\X'_0/\Y$ (this is
the same letter $\L'$ as used for the pullback of the bundle
$\calo(1)$ to the scheme $\X/\Y$). Then $\L'$ is a relatively ample
line bundle on the scheme $\X'/Y$.

\begin{lemma}\label{large.X'.crepant}
The scheme $\X'$ is Cohen-Macaulay and has trivial canonical bundle
$K_{\X'} \cong \calo_{\X'}$. In particular, the map $\pi':\X' \to
\Y$ is crepant.
\end{lemma}

\proof{} We obviously have $K_{\X} \cong \calo_{\X}$, which yields
$K_{\Y} \cong \calo_{\Y}$. To prove that $K_{\X'} \cong
\calo_{\X'}$, it suffices to prove that the canonical rational map 
$$
f = \pi^{-1} \circ \pi':\X' \ratto \X
$$
is defined and injective outside of a subvariety $Z \subset \X'$ of
codimension $\geq 2$. In other words, we have to show that the map
$f:\X' \to \X$ is injective in the generic point of an arbitrary
Weil divisor $E \subset \X'$.

Moreover, since both $\pi$ and $\pi'$ are isomorphisms over the
generic point $\eta \in S$, it suffices to consider divisors $E
\subset \X'$ which lie in the special fiber $X'_o \subset \X'$. 

Let $E$ be such a divisor. Consider the normalization map $\nu:\X'
\to \X'_0$. Since the map $\nu$ is finite, the image $\nu(E) \subset
X'_{0,o} \subset \X'_0$ is also a divisor. But the special fiber
$X'_{0,o} \subset \X'_0$ is irreducible by Lemma~\ref{irr},
therefore such a divisor $\nu(E)$ must coincide with the whole
special fiber $X'_{0,o} \subset \X'_0$. Thus the generic point of
the divisor $E \subset \X'$ lies over the generic point of the
variety $Y$. We are done, since the maps $\pi:X \to Y$, $\pi':X'_o
\to Y$ are by construction both generically one-to-one.

Thus $K_{\X'}$ is trivial. This means that the map $\pi':\X' \to \Y$
is crepant. Since $\Y$ has a smooth crepant resolution $\pi:\X \to
\Y$, the singularities of $\Y$, hence of $\X'$, are canonical. This
implies that these singularities are rational by \cite[II,
Th\'eor\`eme 1]{E}. Hence $\X'$ (and $\Y$) are Cohen-Macaulay.
\endproof

We can now prove that in fact $\X' \cong \X'_o$.

\begin{lemma}\label{fin}
We have
$$
\X' \cong \X'_0 = \Proj_{\Y}\bigoplus_k \pi_*(\L')^{\otimes k}.
$$
\end{lemma}

\proof{} Since the line bundle $\L'$ is relatively ample on
$\X'/\Y$, we have
$$
\X' = \Proj_{\Y}\bigoplus_k \pi'_*(\L')^{\otimes k}.
$$
We now apply Lemma~\ref{emb} to the pair $\pi:\X \to \Y$, $\pi':\X'
\to \Y$ and conclude that
$$
\bigoplus_k \pi'_*(\L')^{\otimes k} \cong \bigoplus_k
\pi_*(\L')^{\otimes k}. \qquad\qquad\qquad\square
$$

Therefore Lemma~\ref{flt} and Lemma~\ref{irr} apply to
$\X'$. Combining them with Lemma~\ref{large.X'.crepant}, we conclude
that the normal irreducible Cohen-Macualay scheme $\X'$ is flat over
$S$ and crepant over $\Y$, and that the special fiber $X'_o \subset
\X'$ is irreducible and reduced.

\subsection{The special fiber.} We have proved all the statements of
Theorem~\ref{main} that deal with the scheme $\X'/\Y$. To finish the
proof, we now have to consider the special fiber $X'_o \subset
\X'$. We have to prove that $X'_o$ is a normal variety crepant over
$Y$, and that
\begin{equation}\label{ravno}
X'_o \cong X' = \Proj_Y\bigoplus_k \pi_*(L')^{\otimes k}.
\end{equation}

\begin{lemma}\label{X'.normal}
The special fiber $X'_o \subset \X'$ is normal.
\end{lemma}

\proof{} Since $\X'$ is Cohen-Macaulay, the hypersurface $X'_o
\subset \X'$ is also Cohen-Macaulay. Thus by Serre's criterion, it
suffices to prove that $X'_o$ is non-singular in codimension
$1$. This splits into two parts.
\begin{enumerate}
\item The scheme $\X'$ is non-singular in codimension $\leq 2$.
\item The function $t$ on $\X'$ which cuts the hypersurface $X'_o
\subset \X'$ has no critical points in codimension $\leq 2$ (in
other words, the differential $dt$ is a non-vanishing $1$-form on a
smooth open subset $U \subset \X'$ whose complement is of
codimension $\codim \X' \setminus U \geq 3$).
\end{enumerate}
To prove these claims, we will use the technique of valuations (see,
e.g., \cite[Section 1]{K2} for an overview). Let $Z' \subset X'_o
\subset \X'$ be an irreducible subvariety of codimension $\codim Z'
= 2$. Consider the valuation $v$ of the fraction field $\C(\X') =
\C(\X) = \C(\Y)$ associated to the exceptional divisor $E$ of the
blow-up of $\X'$ in $Z'$.

It is well-known that the discrepancy of the divisor $E$ in fact
depends only on the variety $\X'$ and on the valuation $v$, not on
the particular geometric realization of this valuation. We will
denote it by $\disc(v,\X')$. Moreover, since both $\X$ and $\X'$ are
crepant over $\Y$, the discrepancy $\disc(v,\X) = \disc(v,\X') =
\disc(v,\Y) = \disc(v)$ is the same when computed with respect to
$\X$, $\X'$ ot $\Y$.

Since $K_{\X'}$ is trivial, the number $\disc(v) = \disc(v,\X') =
\disc(v,\X)$ is either $1$ or $0$. Since the variety $\X$ is smooth,
the center of the valuation $v$ in the variety $\X$ is a closed
subvariety $Z \subset \X$ of codimension $\codim Z \leq \disc(v) + 1
\leq 2$. Moreover, since $\pi(Z) = \pi'(Z') \subset \Y$, the
subvariety $Z \subset \X$ must lie entirely in the special fiber $X
\subset \X$. This excludes the possibility $\codim(Z) = 1$ -- it
would mean $Z = X \subset \X$, hence $\pi'(Z') = \pi(Z) = Y \subset
\Y$ and $Z' = X'_o$, which contradicts $\codim Z'=2$.

Therefore we must have $\codim Z = 2$ and $\disc(v) = 1$. Since $\X$
is smooth, this in turn implies\footnote{This is probably a standard
fact, but an interested reader can find a proof in \cite[Corollary
2.3]{K2}.} that $v$ coincides with the valuation associated to the
exceptional divisor of the blow-up of $\X$ in $Z$. In particular,
since we know that $\X$ is smooth over $S$, this implies that $v(t)
= 1$ for the parameter function $t \in \calo(S)$.

Localizing at the generic point $z' \in Z' \subset \X'$ and using
the theory of algebraic surfaces, we see that $\disc(v)=1$ implies
that $\X'$ is non-singular at $z'$. This proves \thetag{i}. Then
\thetag{ii} immediately follows from $v(t) = 1$.
\endproof

\begin{remark}
Since the scheme $\X'$ is not of finite type over $S$ (nor over
$\C$), the notion of the canonical bundle causes some
problems. Indeed, the sheaf $\Omega^1(\X'/S)$ of relative K\"ahler
differentials is not of finite rank. To make sense of its ``top
exterior power'' $K_{\X'}$, one has to replace $\Omega^1(\X'/S)$
with the quotient sheaf that classifies all derivations compatible
with the $(t)$-adic topology on the local rings of the scheme $\X'$.
This brings no difficulties with the notion of the discrepancy of a
valuation. Indeed, all the valuations that we consider are centered
at the special fiber $X'_o \subset \X'$, hence compatible with the
topology. The same applies to the schemes $\X$, $\Y$.
\end{remark}

\begin{lemma}
The canonical bundle $K_{X'_o}$ of the normal irreducible variety
$X'_o$ is trivial, $K_{X'_o} \cong \calo(X')$. In particular, the
variety $X'_o$ is crepant over $Y$.
\end{lemma}

\proof{} The canonical bundle $K_{\X'}$ is trivial by
Lemma~\ref{large.X'.crepant}. The adjunction formula gives an isomorphism
$$
K_{X'_o} \cong \N(X'_o),
$$
where $\N(X'_o)$ is the normal bundle to the hypersurface $X'_o
\subset \X'$. The differential $dt$ of the parameter $t \subset
\calo(S)$ defines a non-trivial global section of the conormal
bundle to the hypersurface $X'_o \subset \X'$. Moreover, by the
statement \thetag{ii} in the proof of Lemma~\ref{X'.normal}, this
section is non-vanishing outside of a subset $Z \subset X'_o$ of
codimension $\geq 2$. We conclude that the conormal bundle to $X'_o
\subset \X'$ is trivial. Therefore $K_{X'_o} \cong \N(X')$ is also
trivial.  \endproof

\subsection{Identification of the line bundle.}
We see that the variety $X'_o$ is normal and irreducible, and the
projective map $\pi':X'_o \to Y$ is crepant and generically
one-to-one. To prove Theorem~\ref{main}, it remains to establish the
identification \eqref{ravno}.

Denote by $\L'_o$ the restriction of the ample line bundle $\L'$ on
$\X'$ to the special fiber $X'_o \subset \X'$.  We have
$$
X'_o \cong \Proj_Y\bigoplus_{k \geq 0}\pi'_*{\L'_o}^{\otimes k}.
$$
Moreover, we know that $X'_o/Y$ is crepant. Therefore, if we consider
the birational map
$$
f_o = (\pi'_o)^{-1} \circ \pi:X \ratto X'_o,
$$
then Lemma~\ref{emb} shows that
$$
\bigoplus_{k \geq 0}\pi'_*{\L'_o}^{\otimes k} \cong \bigoplus_{k
\geq 0}\pi_*{L_o}^{\otimes k},
$$
where
$$
L_o = f_o^*\L'_o
$$
is the line bundle on $X$ corresponding to the line bundle $\L'_o$
on $X'_o$. Thus to finish the proof of Theorem~\ref{main}, it
suffices to identify the line bundles $L_o$ and $L'$ on the variety
$X$.

To make the statement we have to prove more clear, denote by $i:X
\hookrightarrow \X$, $i':X'_o \hookrightarrow \X'$ the embeddings of
the special fibers. We have a rational map $f:\X \ratto \X'$ and a
rational map of the special fibers $f_o:X \ratto X'_o$. These maps
are compatible with the embeddings: $i' \circ f_o = f \circ i:X
\ratto \X'$. We have an ample line bundle $\L'$ on $\X$. By
definition of the projective scheme $\X'$, the line bundle $f^*\L'$
on $\X$ is the extension to the scheme $\X$ of the given line bundle
$L'$ on the special fiber $X$. In other words, we have $L' \cong
i^*f^*\L'$. On the other hand, by construction we have $L_o =
f_o^*\L'_o = f_o^*(i')^*\L'$. The statement that we have to prove is
$$
f_o^*(i')^*\L' \cong i^*f^*\L'.
$$
This would have been trivial if rational maps $f_o$ and $f$ were
defined everywhere. As things stand, we only have an isomorphism on
the open subset in $X$ where both $f_o$ and $f$ are defined.

Let $Z \subset Y$ be the closed subset of codimension $\geq 3$ that
is provided by Lemma~\ref{codim.3} -- namely, the subset such that
the line bundle $L'$ is relatively generated with respect to the map
$\pi:X \to Y$ over the complement $Y \setminus Z \subset Y$.  By
Proposition~\ref{semismall} the preimage $\pi^{-1}(Z) \subset X$ is
of codimension $\geq 2$.  Therefore to show that $L_o \cong L'$ on
$X$, it suffices to show that $L_o \cong L'$ outside of the subset
$\pi^{-1}(Z)$. Thus it suffices to show that both rational maps $f$,
$f_o$ are defined outside of $\pi^{-1}(Z) \subset X \subset \X$.

The special fiber $\Gamma_o \subset \Gamma$ of the graph $\Gamma
\subset \X \times_\Y \X'$ of the rational map $f:\X \ratto \X'$ can
have several irreducible components. One of these components is the
graph of the rational map $f_o:X \ratto X'_o$. Therefore the
indeterminacy locus of the map $f_o$ is a subset of the
indeterminacy locus of the map $f$. 

We conclude that it suffices to consider the map $f$. Thus to finish
the proof of Theorem~\ref{main}, it suffices to prove that the
rational map $f:\X \ratto \X'$ is defined on the complement to the
closed subset $\pi^{-1}(Z) \subset X \subset \X$. By construction of
the scheme $\X'$, this is equivalent to the following lemma.

\begin{lemma}\label{base}
The base locus $B(\L') = B(\X,\L')$ of the line bundle $\L'$ on $\X$
lies in the closed subset $\pi^{-1}(Z) \subset X \subset \X$.
\end{lemma}

\proof{} We want to show that the image $\pi(B(\L')) \subset \Y$ of
the base locus $B(\L')$ under the projection $\pi:\X \to \Y$ lies in
the closed subset $Z \subset Y \subset \Y$. Since $\Y$ is affine,
the claim is local on $\Y$ -- in other words, the base locus $B(\L')
\subset \X$ coincides with the support of the cokernel of the
canonical sheaf map
$$
\pi^*\pi_*\L' \to \L'.
$$
Denote by $\U = \Y \setminus Z$ the open complement to the closed
subset $Z \subset \Y$, and let $U = Y \setminus Z = \U \cap Y = Y
\setminus Z$ be the complement to $Z$ in the special fiber $Y
\subset \Y$. We have to show that the canonical map $\pi^*\pi_*\L'
\to \L'$ is surjective on the preimage $\pi^{-1}(\U) \subset
\X$. Since the projection $\pi:\X \to \Y$ is one-to-one outside of
the special fiber $X \subset \X$, the map $\pi^*\pi_*\L' \to \L'$ is
tautologically surjective outside of the special fiber. Therefore it
suffices to prove that this map is surjective on the special fiber
$\pi^{-1}(U) = \pi^{-1}(\U) \cap X \subset \pi^{-1}(\U)$.  More
precisely, we have to prove that the map
$$
i^*\pi^*\pi_*\L' \to i^*\L'
$$
is surjective on the subset $\pi^{-1}(U) \subset Y$, where $i:X \to
\X$, $i:Y \to \Y$ denotes the embedding of the special fiber.

By definition we have $i^*\L' = L'$, and we tautologically have
$i^*\pi^*\pi_*\L' \cong \pi^*i^*\pi_*\L'$. But the canonical bundle
$K_X$ is trivial, and the line bundle $L'$ is by assumption
relatively generated over $U \subset Y$. By Grauert-Riemenschneider
vanishing this implies that
$$
R^k\pi_*L' = R^k\pi_*K_X \otimes L' = 0
$$ 
on $U \subset Y$ for all $k \geq 1$. Then the base change gives a
line bundle isomorphism
$$
\pi^*i^*\pi_*\L' \cong \pi^*\pi_*i^*\L' \cong \pi^*\pi_*L',
$$
and the map $\pi^*\pi_*L' \to L'$ is by assumption surjective on the
open subset $\pi^{-1}(U) \subset X$.
\endproof

\begin{remark}
Lemma~\ref{base} is the first place in our construction where we
have used the assumptions on the line bundles $L'$ other than
Condition~\ref{techn}. In fact, we can start with an arbitrary line
bundle $L'$. Then everything up to Lemma~\ref{base} works just as
well. However, Lemma~\ref{base} does break down. If the original
line bundle does not satisfy some form of Grauert-Riemenschneider
vanishing, then the indeterminacy locus of the rational map $f:\X
\ratto \X'$ may contain a divisor $E \subset X$ in the special fiber
$X \subset \X$. In this case the graph $\Gamma \subset \X
\times_{\Y} \X'$ of the map $f$ has at least two irreducible
components of dimension $\dim X$ in the special fiber $\Gamma_o
\subset \Gamma$: one is the graph of the rational map $f_o:X \ratto
\X'_o$, and the other lies entirely over the divisor $E \subset X$.
\end{remark}

\section{An application -- smoothness of flops}\label{app}

We finish the paper with one concrete geometric corollary of
Theorem~\ref{main} -- both because it is interesting in its own
right, and in order to convince the reader that even in spite of the
very restrictive Condition~\ref{techn}, our result does have
immediate applications. This application is
Proposition~\ref{flp.smth} -- the smoothness of flops. To save
space, we only give a sketch of the proof (omitting exact references
for deformation theory statements and the like).

\subsection{Smoothness of flops.}\label{flp.smth.sub}
We note that Proposition~\ref{flp.smth}~\thetag{i} is an immediate
corollary of Lemma~\ref{pic.rigid}; the real problem is to prove
\thetag{ii}. We will actually prove a slightly more general
statement. For this we return to the setting of Theorem~\ref{main}.

\begin{prop}\label{smth}
In the setting of Theorem~\ref{main}, assume that the scheme $\X'$
is $\Q$-factorial. Then the scheme $\X'$ and the special fiber $X'
\subset \X'$ are smooth, and the special fiber $X' \subset \X'$ is
symplectic.
\end{prop}

\proof{} By definition of the twistor deformation
(Definition~\ref{tw.def}), the scheme $\X$ carries an ample line
bundle $\L$ such that
$$
c_1(\L) = \Omega \cntrct \theta_{\X},
$$
where $\Omega \in \Omega^2(\X/S)$ is the relative symplectic form on
the twistor deformation $\X/S$, $\theta_X \in H^1(\X,\T(\X/S))$ is
the Kodaira-Spencer class of the deformation $\X/S$, and $c_1(L) \in
H^1(\X,\Omega^1(\X/S))$ is the relative first Chern class of the
line bundle $\L$.

The rational map $f:\X' \ratto \X$ is defined on an open subset $U
\subset \X'$ whose closed complement $\X' \setminus U \subset \X'$
is of codimension $\geq 2$. Moreover, as the proof of
Lemma~\ref{irr} shows, we can actually assume that $f$ is injective
on $U \subset \X'$. In particular, $U$ is smooth over $S$, and we
have
$$
c_1(f^*\L) = f^*\Omega \cntrct \theta_U,
$$
where $\theta_U \in H^1(U,\T(U/S))$ is the Kodaira-Spencer class of
the deformation $U/S$. Since $\X'$ is by assumption $\Q$-factorial,
some multiple of the line bundle $f^*\L$ extends to the whole
$\X'$. Therefore the Chern class $c_1(f^*\L) \in
H^1(U,\Omega^1(U/S))$ extends to a class in
$H^1(\X',\Omega^1(\X'/S))$.

The maps $\pi:\X \to \Y$, $\pi':X' \to \Y$ are one-to-one over the
generic point $\eta \in S$. Therefore it suffices to prove that the
schemes $\X'$ and $X'$ are regular in every closed point $x \in X'$.

Let $x \in X'$ be such a point, and let $V \subset \X$ be an open
affine neighborhood of the point $x \in X' \subset \X'$. Since $V$
is affine, every class in the cohomology group
$H^1(\X',\Omega^1(\X'/S))$ vanishes after restriction to $V \subset
\X'$. This applies in particular to the class
$c_1(f^*\L)$\footnote{The scheme $V$ is not {\em a priori} smooth
over $S$. However, one can still define the Chern class $c_1(f^*L)
\in \Omega^1(\X'/S)$ -- for instance, by applying $d\log$ to the
corresponding class in $\Pic(V) = H^1(V,\calo_V^*)$.}. We conclude
that $c_1(f^*\L) = 0$ on the intersection $U \cap V$. Since the
relative symplectic form $f^*\Omega \in \Omega^2(U/S)$ is
non-degenerate, this implies that the Kodaira-Spencer class
$$
\theta_{U \cap V} \in H^1(U \cap V,\T(U/S)) = \Ext^1_{U \cap
V}(\Omega^1(U/S),\calo(U \cap V))
$$
also vanishes.

But the scheme $\X'$ is normal, and the complement to the subset $U
\cap V \subset V$ is of codimension $\geq 2$. This means that if we
denote by $j:U \cap V \hookrightarrow V$ the embedding, then the
direct image $j_*\calo(U \cap V)$ coincides with the structure
sheaf,
$$
j_*\calo(U \cap V) \cong \calo(V). 
$$
Therefore for every complex $\F$ of coherent sheaves on $V$ which is
concentrated in non-positive degrees we have
\begin{align*}
\Ext^1_V(\F,\calo(V)) 
&= \Ext^1_V(\F,j_*\calo(U \cap V)) \subset \\
&\subset \Ext^1_V(\F,R^\hdot j_*\calo(U \cap V)) 
= \Ext^1_{U \cap V}(j^*\F,\calo(U \cap V)).
\end{align*}
Applying this to the relative cotangent complex $\Omega^\hdot(V/S)$
of the flat deformation $V/S$, we conclude that the Kodaira-Spencer
class
$$
\theta_V = \theta_{U \cap V} \in \Ext^1(\Omega^\hdot(V/S),\calo(V))
\subset \Ext^1_{U \cap V}(\Omega^\hdot(V/S),\calo(V \cap U))
$$
of this deformation is trivial. Therefore for every integer $n \geq
0$, the associated order-$n$ infinitesimal deformation $V_n$ of the
special fiber $V \cap X'$ is trivial. Moreover, we can choose a
compatible system of projections $V_n \to V_o$ onto the special
fiber $V_o = V \cap X' \subset \X'$.

Passing to the limit, we obtain a projection $\tau:V \to V_o$, hence
a product decomposition $V = V_o \times S$. Let $\wt{x}:S \to x
\times S \subset V_o \times S$ be the horizontal section of the map
$\X' \to S$ which passes throught the point $x \subset X'$. Then the
conormal sheaf $\N(\wt{x})$ to this section is a constant vector
bundle on $S$. Its rank is equal to $\dim T_xX'$, the dimension of
the Zariski tangent space $T_xX'$ to $X'$ at $x$.

But since $V$ is by assumption smooth over the generic point $\eta
\in S$, the rank of the normal bundle $\N(\wt{x})$ must concide
with $\codim(\wt{x},\X') = \dim(X')$, and we have $\dim T_xX' =
\dim(X')$. This means that the scheme $X'$ is regular at the point
$x \subset X'$. Since the point $x \subset X'$ was arbitrary, we
conclude that both $\X'$ and $X'$ are smooth algebraic varieties.

Finally, we note that for every $k \geq 0$ the spaces of $k$-forms
are the same for all smooth birational models of the same algebraic
variety $X$ (it suffices to prove this for blow-ups with smooth
centers, where it immediately follows from a direct
computation). Therefore the symplectic form $f^*\Omega$ extends from
the open subset $U \cap X' \subset X'$ to the whole $X'$.
\endproof

\proof[Proof of Proposition~\ref{flp.smth}.] By the uniqueness of
flops (\cite[Proposition 5-1-11(2)]{KMM}), we can assume that the
flop $\X'$ of the twistor deformation $\X/\Y$ is the same scheme as
the one provided by Theorem~\ref{main} (Condition~\ref{techn} is
satisfied precisely because the flop $\X'$ is assumed to
exist). Then the scheme $\X'$ is $\Q$-factorial by \cite[Proposition
5-1-11(1)]{KMM}, so that we can apply Proposition~\ref{smth}.
\endproof

\bigskip

\noindent
{\em E-mail address:\/} {\tt kaledin$@$mccme.ru}

\end{document}